\documentclass[12pt,a4paper,reqno,dvipsnames]{amsart}
\usepackage{tcolorbox}
\usepackage{xcolor}
\usepackage[hyphens]{url} 
\usepackage{amsmath}
\usepackage{mathtools}
\usepackage[utf8]{inputenc}
\usepackage[T1]{fontenc}
\usepackage[english,greek]{babel}
\usepackage{textgreek}
\usepackage{CJKutf8}
\usepackage{accessibility}
\usepackage{ulem}

\makeindex

\begin{document}
\selectlanguage{english}

\title[Human Thought in the Age of AI]{Mathematical methods and human thought in the age of AI}

\author{Tanya Klowden and Terence Tao}

\begin{abstract}
    Artificial intelligence (AI) is the name popularly given to a broad spectrum of computer tools designed to perform increasingly complex cognitive tasks, including many that used to solely be the province of humans. As these tools become exponentially sophisticated and pervasive, the justifications for their rapid development and integration into society are frequently called into question, particularly as they consume finite resources and pose existential risks to the livelihoods of those skilled individuals they appear to replace.
    In this paper, we consider the rapidly evolving impact of AI to the traditional questions of philosophy 
    with an emphasis on its application in mathematics and on the broader real-world outcomes of its more general use. We assert that artificial intelligence is a natural evolution of human tools developed throughout history to facilitate the creation, organization, and dissemination of ideas, and argue that it is paramount that the development and application of AI remain fundamentally human-centered. With an eye toward innovating solutions to meet human needs, enhancing the human quality of life and expanding the capacity for human thought and understanding, we propose a pathway to integrating AI into our most challenging and intellectually rigorous fields to the benefit of all humankind.
\end{abstract}

\maketitle

\section{Introduction}

It is a testament to how quickly artificial intelligence (AI) technologies have been deployed into every crevice of digital life that, in the process of composing this paper using standard tools, the authors had no less than three different digital agents insert themselves into the narrative unsolicited\footnote{Any and all of these AI ``contributions'' were promptly removed from the text.}. Humanity is standing at the threshold of a digital Industrial Revolution, unfolding at unprecedented speeds. In the physical sciences, AI advances have led to Nobel prize-winning research \cite{alphafold}; while in the humanities, fears abound that the generative text capabilities of modern AI could be the death of the subject\cite{marche}. As language translators have thrown the doors wide open for cultural exchange and international cooperation, a flood of deepfakes and slop has followed, sloshing through our digital third spaces. AI quickly went from a novelty, to a vital resource, to (in some cases) a present existential threat \cite{anguiano}.

\subsection{Our definition of artificial intelligence}

For the purposes of this article, AI refers to the broad spectrum of computer tools designed to perform increasingly complex cognitive tasks, including many that used to solely be the province of humans.  AI tools are extremely diverse, ranging from the data-driven \emph{machine learning} (ML) technologies of today (such as \emph{large language models} (LLMs) that can process complex text, or \emph{diffusion models} that can generate images and other media), to the more traditional \emph{good-old fashioned AI} (GOFAI) (such as automated theorem provers or chess engines), which can solve narrow ranges of problems by applying precise mathematical rules.

\subsection{Purpose of this article}
There has been no shortage of discussion about what these tools can or cannot do; but comparatively less discussion of \emph{why} these tools are being so rapidly developed and deployed or \emph{how} they impact the billions of lives that interact with them for research and education, for work, for play, and even for rest\cite{oh}. The authors of this paper come from academic domains that are frequently viewed as polar opposites: mathematics and the study of art.  But we both have found it beneficial to incorporate several AI tools into our disparate areas of research on a day-to-day basis, and found a surprising amount of common ground regarding the very messy, but universal, philosophical questions that real-world AI use poses. Using mathematics as a model, we will consider the benefits, risks, ethics and outcomes of incorporating AI into routine workflows and then expand these observations to broader real-world use. Despite the risks that these new, and not necessarily morally neutral technologies present, we argue two-fold that AI tools should be developed, implemented, and applied both within mathematics and in other domains: they have the potential to radically augment our natural human abilities and they are capable of expanding what is possible beyond what we humans could do individually or within the limits of our own collective capacity. Drawing from our own experiences with these tools, we particularly examine the human/AI interface and offer suggestions on the evolution of these technologies in ways that offer more benefits than harms to humanity and value the unique contributions of human thought and action in concert with the new modalities that future AI development promises. 

\subsection{The Faustian Bargain}

The incentives of market competition have fueled a frenzied pace of development of AI technologies and fascinated entire industries with visions of radically accelerated workflows and cost savings.  The ``prisoner's dilemma'' of such competition has pressured many individuals and organizations to experimentally adopt these tools as hastily as possible, at the expense of more deliberate evaluation of the economic, social, or moral costs and benefits of such an adoption -- or, more fundamentally, why we should be developing such technologies in the first place.  As such, we collectively have already adopted \emph{de facto}  a ``Faustian bargain'' with these technologies, giving them increasing access to our data, cognitive workflows, and decision processes, in exchange for the promise of being able to accomplish a greater range of tasks with increasing efficiency and with less tedious effort.  

In theory, technology is morally neutral; it can empower both positive and negative use cases.  But through this empowerment, it exacerbates existing moral dilemmas, and creates new ones.  For instance, the horrific medical research on prisoners during World War II which led to lifesaving data on the limits of human endurance, raised difficult questions regarding the ethicality of using such data to develop new medical advances \cite{swain} .  While not as gruesome, the murky provenance of the data and intellectual property used to train the current generation of AI tools arguably raises similar questions today \cite{tarkowski}.

When a technology develops slowly enough, it is possible to have the necessary philosophical conversations and debates about it before it is widely deployed; stem cell research is one notable example of this.  However, modern AI technologies are already widely deployed, with no practical way to ``put the genie back in the bottle''; ironically, strict regulation imposed at this point would disproportionately shut down the more positive use cases of AI, such as in the acceleration of scientific research, without eliminating the more wasteful or malicious uses of the technology.  Pragmatically, the discussion about AI has now moved towards how to manage coexistence with these technologies: evaluating the costs and benefits of AI (both in academic disciplines, and in to broader society), and identifying best practices and frameworks to use AI in as positive a way as possible, while simultaneously discouraging the (many) ways in which these tools can be used poorly to degrade the reliability and value of our cognitive achievements.

\section{Historical parallels: is this time different?}

\subsection{Past integration of automation technologies}

Automation is of course not a new phenomenon.  Many past technologies have also enabled the ability to automate tasks previously assigned to humans, eliminating or greatly reducing the need for some types of human jobs, while creating or increasing the need for others; in some cases.  Within the scientific community for instance, ``phase transitions'' have occurred in which broadly and rapidly switched over to new tools (such as the internet, the use of computers for scientific computation, or even humble typesetting languages such as \LaTeX) due to their evident advantages.  
But these past technologies have mostly affected secondary aspects of the profession, such as the communication and dissemination of results rather than the creation of such results.  
And, while the tasks automated by these tools required specialized training and expertise to perform, they typically did not require an understanding of more philosophical aspects of a profession, such as the nature of knowledge, beauty, meaning, or morality  \cite{chun}.  Of course, such technologies could still generate  discussions on philosophical topics -- for example, whether there were inherent aesthetic or creative features of an original piece of art that no mechanical reproduction could properly capture, or on the moral and ethical implications of the displacement of labor caused by the Industrial Revolution -- but they were not considered to contest the fundamental philosophical assumptions underlying such discussions.  For instance, the invention of the printing press revolutionized the communication of information and ideas, but it did not significantly alter the understanding of what an idea or a piece of information \emph{was}; the original generation of this content was still performed by the deliberate actions of humans.

\subsection{Modern AI}

But modern AI can automate large portions of the creative process itself, enabling the mass-generation of intellectual products, such as artwork, mathematical proofs, or scientific or philosophical theories, with far less human oversight than was previously required\footnote{Current tools typically still require a human to generate an initial prompt for the AI to follow, but this process can itself now be largely automated as well.}.  This has created an unprecedented decoupling between the outward \emph{form} of such products, and the values and thought processes used to \emph{create} these products.  A diffusion model may now create an aesthetically pleasing landscape, for instance, which was not directly inspired by any particular location in the physical world, though countless images of actual landscapes (as well as many images completely unrelated to landscapes) were certainly used to train the outputs of that model; the aesthetic response of the image thus becomes decoupled from the original sources of such aesthetics.

This is not new philosophical territory by any means. Searle's ``Chinese room'' thought-experiment \cite{searle}, regarding the question of whether a mechanical device programmed to converse in Chinese truly understands the language, dates back to 1980.  The ``AI effect'' also was recognized around this time; for instance, the ability to perform well at chess was considered a good measure of intelligence until the advent of chess engines which could ``mindlessly'' outperform chess masters through mechanical exploration of game trees, at which point the ``chess test'' for intelligence became largely abandoned.  The famed ``Turing test'' of whether an AI could converse in a manner indistinguishable from humans has similarly been effectively passed by modern LLMs (see, e.g., \cite{turing}), relinquishing its former status as a ``gold standard'' for artificial intelligence.   For a more recent discussion, see \cite{Chen}.

For now, we can still point to markers of ``fundamental'' understanding, such as the ability (or lack thereof) to coherently explain and defend the creative processes that led to a new artwork, mathematical proof, or other intellectual product, as a still-viable test to distinguish between human and AI-generated content, but if future generations of AI somehow also manage to convincingly pass such tests as well, would we have to move the goalposts once again on what intelligence, understanding, and creativity actually \emph{are}?  Would the definitions, values, and objectives of such disciplines as mathematics and the humanities need to be re-evaluated?  And what status should we grant these increasingly sophisticated AI tools - will they be assistants, co-authors, or even independent creators in their own right?  And if so, how should we treat the content they produce, and the intellectual processes that led to such content?

\section{Mathematics as a sandbox for AI use}

Such broader philosophical questions about AI are extremely complex and multifaceted, and we of course do not pretend to have definitive resolutions to any of them; and the speed of change in this space is such that any proclamations we make are at risk of being overtaken by striking new technological advances.  However, we can offer some perspectives from the world of mathematics, both in the realm of pure mathematical reasoning, and in the emerging application of modern mathematical analysis in the humanities.
We view mathematics as a suitable ``sandbox'' for exploring broad questions such as the impact of AI across the sciences (and society as a whole), as it has an older and more advanced foundation, and is by its nature well suited to explore a variety of hypothetical abstract scenarios which are counterfactual to reality.  It is our hope that the lessons learned from integrating (or not integrating) AI into mathematics can give broader perspectives on how AI will interact with sciences and society in general.

Frontier AI models can now solve increasingly complicated mathematical problems, with proofs that can be independently verified, without directly reproducing the problem-solving practices of human mathematicians (such as testing out special cases, and then generalizing from those examples), though its training data would include proofs generated in such a traditional fashion; and so mathematicians will increasingly encounter situations in which the ability to prove theorems is decoupled from the reasoning processes needed to discover and understand such proofs.  This contributes to an existing trend of decentralization in modern mathematics; in a world where advanced mathematics is needed in an extremely broad range of applications, the ``Bourbaki era'' \cite{mashaal} of having a central authority prescribe the orthodox practice of mathematics is already decades in the past\footnote{Though one could argue that the ongoing projects to create large unified libraries of formal mathematics, such as Lean's \emph{Mathlib} project, could be a modern successor to the efforts of the Bourbaki group.}.

At the current state of the technology, the most sophisticated AI tools still exhibit significant and often bizarre weaknesses; they can achieve remarkable and super-human performance in some tasks, while simultaneously demonstrating often hilarious levels of basic misunderstanding and error in others.  Mathematics is no exception to this phenomenon. AI-generated mathematics can appear superficially flawless - which is to be expected, since these models are designed to produce outputs as visually close to correct human-generated proofs as possible - while also making fundamental mistakes (for instance, asserting that all odd numbers are prime) that would have been trained out of a human mathematician at an early stage of their training, and can often make the resulting argument unsalvageably nonsensical.  At the same time, this top-down approach of focusing primarily on generating good-looking outputs rather than on the fundamental cognitive processes that were traditionally needed to create such outputs can be surprisingly effective; the same AI that routinely makes basic mathematical errors, can also mysteriously arrive at the correct answer to a complex math problem with superior accuracy to human experts, or even supply a strange but technically correct proof that the answer is valid.

Significant effort is now being directed to reduce or eliminate these weaknesses of AI as much as possible; often not by directly strengthening the AI's innate ``understanding'' of any given intellectual task, but rather by placing such AI tools in a rigorous environment of independent testing, training, and verification to reduce the numerical incidence of errors.  The ability to resolve deep mathematical conjectures
is still currently out of reach of a completely autonomous AI, but it is very plausible in the near future that such AI tools could greatly assist human mathematicians in such endeavors, even if we would still hesitate at describing such assistance as the expression of genuinely intelligent thought.
Still, the fact that such mechanistic and error-prone approaches to as intellectual a discipline as mathematics can (or soon will) generate so many of the traditional markers of quality in the subject indicates that we have to re-evaluate our models of what intelligence or creativity actually \emph{is}, and how it is to be measured.

\section{AI and the nature of mathematical truth}

\subsection{Mathematics and standards of proof}

Mathematics\footnote{Here we leave vague the concept of what mathematics actually \emph{is}.  One can adopt a prescriptivist view, for instance using the Davis--Hersh definition \cite{davis} of mathematics as ``the study of mental objects with reproducible properties''.  Or one could adopt a descriptivist view, namely that mathematics is the activity that mathematicians actually perform in practice.  Our discussion here somewhat favors the latter view.} has had a long tradition of having an objective standard of proof, starting with Euclid and refined with the advent of stable and (empirically) secure foundations of mathematics in the early twentieth century.  It has been noted (see, e.g., \cite{consensus}) that the near-universal acceptance of these foundations has given modern mathematics the rare and precious ability to arrive at a consensus on the validity of any given argument or assertion in the field, since (in principle) one could insist that such arguments be spelled out in such fine detail that each individual step could be checked to be a correct application of the standard axioms and logical inference rules of mathematics.  A typical example was the claim by Nelson \cite{nelson} in 2011 that the Peano Axioms were logically inconsistent; this was a claim very far from the mathematical mainstream, and yet it was possible to resolve the issue by pointing out a subtle flaw in the argument that Nelson readily accepted, thus withdrawing the claim.

However, in practice, the arguments of human mathematicians fall short of the ideal of perfectly rigorous proof; minor and major mistakes in the literature are common, with some of these being corrected by formal errata or revisions, and others being neglected or passed on informally as part of the ``folklore'' of the subfield.  Arguments which are heuristically plausible are often accepted with minimal checking, while surprising assertions that go against the conventional wisdom are met with heavy skepticism, even if the arguments do ultimately turn out to be correct on a line-by-line reading.

\subsection{The Smell Test}

Until now, this state of affairs has been reasonably satisfactory; human mathematicians who follow good heuristics and intuition tend to produce convincing proofs that are largely correct, with most errors fixable, whereas mathematicians who lack such intuition tend to produce proofs that contain enough superficial issues that one can be rightfully suspicious of the content even before one checks carefully.  Informally, human-generated mathematical arguments tend to come with\footnote{We use this sensory metaphor in analogy with the concept of a ``code smell'' in software engineering.} a ``smell'' that experienced mathematicians (perhaps subconsciously) use obtain their initial impressions of how convincing the argument is, well before they have been able to check the individual steps of that argument.  For instance, the blog post ``Ten signs a claimed mathematical breakthrough is wrong'' of Aaronson \cite{aaronson} lists some common examples of arguments that exhibit such a ``bad smell'', that one can detect well before one has located a specific logical flaw in the proposed argument. And not all errors are equally disastrous; some errors may even have some beneficial value, for instance by revealing a promising approach before being able to fully validate it \cite{deo-2}.

One component of a favorable ``smell'', as noted\footnote{See also the article \cite{tao} by the second author, which argues that ``good'' mathematics, regardless of how it is initially defined, often tends in practice to fit into broader mathematical narratives, such as the dichotomy between structure and randomness, or the ability of algebra to explore questions about geometry (or vice versa).} by Thurston \cite{thurston}, is the sense that an argument is providing \emph{understanding} or \emph{insight}; that it not just shows that a certain set of hypotheses logically entails a given conclusion, but provides a causal narrative for how that entailment was possible, and which parts of the argument were performing the ``heavy lifting'', which parts were novel or surprising compared to previous literature, and which ones were routine technical considerations.  Such interpretations and impressions of the mathematical text generally are not captured in the official frameworks of rigorous mathematics, such as first-order logic or set theory; but they are essential in allowing the human mathematicians reading the argument to draw broader lessons about how one would expect the arguments to generalize to other settings, or interact with other methods in the literature.  Such narrative structures also help strengthen confidence in the robustness of an argument; a single misplaced sign in a calculation could invalidate a lengthy mathematical argument, but if the proof had a clear strategy regarding how the key difficulties in the argument were systematically isolated and addressed, following analogies with previous successful arguments in the literature, then it becomes more likely that local errors in the argument can be repaired while staying true to the spirit of the original proof.

\subsection{Formalization to the rescue?}

There are several developments that may force the mathematical community to re-evaluate this semi-formal standard of proof.  One of them is of a technical nature: as mathematics has matured and become more sophisticated (and increasingly computer-assisted), the arguments have become longer and more complex, with cutting-edge papers in some fields routinely exceeding a hundred pages in length, making line-by-line verification by human referees increasingly onerous.  In practice, this has meant that such careful checking is not always performed, save for the most high-profile and important results, leading to an increasing (over-)reliance on the aforementioned sense of ``smell'' to assess the credibility of mathematical arguments.

It seems possible that such issues could be resolved (or at least ameliorated) by technical means, and in particular through the more widespread deployment of \emph{formal proof assistants} (such as \emph{Lean} or \emph{Rocq}) which can automatically check the validity of a mathematical argument if it is written in a certain precise computer language \cite{toffoli}.  Such formalization remains too tedious at present to deploy systematically (converting a traditional, informally written proof into such a formal language typically takes about five to ten times longer than writing that proof in the first place), but there are significant efforts underway to make the process faster and more user-friendly, for instance by integrating AI tools to achieve partial (or possibly even complete) ``autoformalization'' \cite{wu}.

\subsection{Limitations of formal verification}

But even if such technical issues are resolved, and mathematical proofs routinely come up with a formal verification of correctness, several new issues arise, especially in a near future in which increasingly sophisticated arguments may be partially or fully generated by AI tools.  Firstly, formal verification only certifies that a formalized argument establishes a formal mathematical statement, but does not rule out errors in translation between the formal statement and the original intended statement.  For instance, Fermat's last theorem asserts that for any natural number $n$ greater than $2$, there are no natural number solutions $a,b,c$ to the equation $a^n + b^n = c^n$; but implicit in this informal description is the convention that the natural numbers start at $1$ rather than $0$.  An AI tasked to solve this problem may erroneously assume that $a,b,c$ are permitted to be zero, and based on this produce a (formally certified) proof that Fermat's last theorem is false!  Thus, while formalization can in principle significantly reduce the need to perform careful human review of informal mathematical text, it does not eliminate the need for such review entirely\footnote{It is even theoretically possible that mathematics itself could be ``hacked'' by subtly manipulating the formalization of key definitions in standard formal mathematical libraries such as \emph{Mathlib}; see \cite{tanswell}.  Ironically, the increasingly collaborative, social, and large-scale of mathematical research, while generally a highly positive development, may also increase potential vulnerabilities to such attacks that were not a significant concern in prior eras when mathematics was mostly performed by small groups of individuals.}.

Secondly, even in the purely abstract setting of advanced mathematics, only a portion of a given argument can be formulated in the type of deductive logic that is amenable to formalization.  While deductive proof remains the crucial core of most mathematical work, there is a penumbra of heuristic, empirical, or metamathematical reasoning around this core which provides valuable information on \emph{why} the argument works, \emph{whether} it extends to other contexts, \emph{what} the motivation is for pursuing these questions, and \emph{how} one might reconstruct the argument from more basic principles.  Human-written proofs, by their nature, will tend to provide this penumbra organically as part of the writing process (particularly of the authors are skilled at exposition); but an AI that has been trained specifically on the criterion of formal correctness, at the expense of all other considerations, could produce ``odorless'' proofs which superficially resemble a well-written human proof, and may even pass formal verification tests, but yet remain strangely unsatisfying - fulfilling to the letter the explicit objective to establish a given mathematical claim, while yielding far less insight than expected on the broader mathematical field that this claim is part of.  In a world where all media generated is AI-polished to a high sheen, including mathematical proofs with beautiful typesetting and clear, GPT-produced explanations, is something lost in forsaking the grubbier, messier world of hand-written (or at least hand-typed) text?

\subsection{Adaptation to earlier challenges}

The mathematical community has adapted to previous technological challenges to its standard of proof.  Large computer-assisted proofs, such as the proof of the four color theorem \cite{appel} or the Kepler conjecture \cite{kepler}, were initially quite controversial, being impractical to fully check by hand; but in time new standards of establishing confidence were established for these types of arguments, such as providing replicable code, isolating the computational components of an argument in specific, clearly stated lemmas separate from the more conceptual aspects of a paper, and providing additional related data and ``checksums'' to check that the computer-generated calculations agree with various ``sanity checks''.  In effect, these developments shifted the standards of proof in mathematics in the direction of that in the natural sciences, in which both theoretical argument and empirical experiment, when properly designed, executed, and reported, are acceptable sources of scientific truth.

\subsection{The evolution of AI-assisted mathematics}

Similar evolutions will take place\footnote{Our thinking here has been influenced by the views of other mathematicians on this topic, including \cite{venkatesh}, \cite{deo}, \cite{deo-2}, \cite{avigad}, as well as the broader discussion in \cite{ams-1}, \cite{ams-2}.} with the advent of significant AI-assisted or AI-generated mathematics.  The burden of producing verified deductive proofs may increasingly fall to computers rather than humans, with proofs increasingly being restructured\footnote{For several concrete examples of such restructuring and further exploration of these developments, see \cite{macbeth}.} so that tedious calculations that would previously be carefully arranged to be human-verifiable are increasingly outsourced to software tools instead. For instance, infamous phrases in mathematics such as ``the proof is left to the reader'' or ``By standard arguments, we have'', for instance, would instead be replaced with a call to an LLM that produces both human-readable and computer-verifiable justifications for such claims.
With advances in auto-formalization, it will also become significantly easier to explore how a given argument changes with respect to specific choices of foundations of mathematics, allowing for the metamathematics\footnote{One example of such metamathematics is the \emph{reverse mathematics}  (see, e.g, \cite{stillwell}) of a theorem, which seeks to understand which axioms of mathematics (e.g., the axiom of choice, or the law of the excluded middle) are actually needed to establish a given result.  Traditionally, the reverse mathematics of a result is only explored many years after the original proof of the result, and requires specialist training in logic as well as domain expertise for the subfield of mathematics that the theorem resides in.} of a result to be rigorously discussed and explored simultaneously with the mathematical result itself.

At the same time, more focus and attention may be given in the future by human mathematicians to ``softer'' aspects of mathematical reasoning, such as heuristics and motivation for pursing a result or selecting a proof strategy for that result, experimental evidence\footnote{In particular, given the increasing ability of AI to be able to ``guess'' the answer to even extremely complex mathematical questions without having anything resembling a formal proof, it will become increasingly necessary to develop standard procedures for citing and incorporating such unverified guesses into the mathematical literature in a responsible fashion.} in favor of (or against) the result, or the trial-and-error process leading to the discovery of a working argument.  These aspects are not as easy to automatically verify and measure as deductive proof, and thus less amenable\footnote{Another hurdle to automating these aspects of the mathematical research process is a relative lack of data; published literature tends to focus on successful proofs of results, at the expense of detailing the (often quite rich and nuanced) processes, both formal and informal, that led to such proofs.} to machine learning strategies such as reinforcement learning.  It is conceivable that professional mathematicians may increasingly adopt\footnote{In particular, one can envision an increasing division of labor in the future of mathematical research: while all mathematicians should stay broadly familiar with the different stages of proposing, establishing, and then interpreting mathematical results, any given mathematician may increasingly specialize in just a few aspects of this process, for instance focusing on utilizing AI assistants to prove results as directed by some more senior member of a research group, or on using the most recent literature produced by some combination of human mathematicians and AI assistants to propose new directions of inquiry.} modes of argument from other disciplines, such as the theoretical and experimental sciences or even the humanities, to buttress their core deductive arguments with additional types of reasoning, such as statistical analysis of experimental data, or speculative theorizing guided by both confirmed mathematical results and non-rigorous philosophical principles.  Historically\footnote{For instance, a previous proposal by Jaffe and Quinn \cite{jaffe-quinn} to systematically develop a field of ``theoretical mathematics'' received a largely negative reception from professional mathematicians, leading to multiple rejoinders including the aforementioned article of Thurston \cite{thurston}.}, mathematicians have been reluctant to stray too far from their ``gold standard'' of rigorous deductive proof, due in part to the many visible examples of low-quality mathematics that can be produced when one no longer adheres to such standards\footnote{Kim \cite{kim} invokes a currency metaphor to describe the social dynamics: professional mathematicians need to accumulate some credibility ``currency'', by proving difficult new mathematical results, before they can ``afford'' to ``spend'' that currency on speculative activities, such as formulating conjectures or philosophizing about the broader consequences of a result.}; but in a future era where proofs can be automatically generated and verified in a highly trusted fashion, there may be more opportunity to safely explore such broader modes of mathematical reasoning and discussion.

These new technologies could also impact the longer-term goals of mathematics in significant negative ways.  At the educational level, we are already seeing many students who resort almost immediately to modern AI tools to perform their assigned coursework, achieving the immediate goal of producing verifiable answers to a given problem at the expense of developing more sustainable mathematical skills and intuition; similarly at the research level, the ``fourth paradigm'' of data-driven mathematics \cite{fourth} could conceivably be so successful as to crowd out the more traditional paradigms of empirical evidence, theoretical reasoning, and computational numerics (the second of which being the currently dominant paradigm for pure mathematics), as well as the great value that human mathematicians\footnote{Somewhat related to this, aesthetic notions such as the ``beauty'' or ``elegance'' of a mathematical argument may become even more decoupled than they currently are from the formal correctness of such arguments. Consider for instance the proofs generated by \emph{AlphaProof} \cite{alpha} to problems in the 2024 International Mathematical Olympiad, which contained numerous redundant or inexplicable steps but nevertheless were formally verified to be correct solutions.  See also the discussion in \cite{deo}.} gain from visual, kinesthetic, and other sensory intuition, or from intuition grounded by our familiarity with the laws of physics, economics, biology, etc..  Even assuming a completely trusted implementation of formal methods, an uncritical embrace of AI assistance in the mathematical research space could lead to the undesirable outcome of a flood\footnote{This could be viewed as an illustration of the law of unintended consequences.  In past mathematical eras where the task of obtaining a rigorous mathematical proof required painstaking human effort, mathematical activity naturally focused on problems which were deemed by the mathematical community to be of interest, even if the philosophical question of what it truly meant for a given  result to be ``interesting'' or ``relevant'' was often not explicitly considered by most members of that community; the evolution of the literature was slow enough that this largely social mechanism of determining mathematical significance could correct itself over time.  In a future era where mathematical results can be mass-produced at significantly faster speeds due to automation, such philosophical issues may require far more active attention.  See also \cite{rittberg}, \cite{trust} on the need to make value judgments, including trustworthiness of the author, when deciding whether to allocate attention a claimed mathematical result.} of largely AI-generated papers containing results that are technically correct and new, but which do not contribute to broader mathematical narratives, and do not build up intuition for either the authors or the readers.  The negative impressions produced by such low-quality work may lead to a stigma against even the most careful and responsible application of AI assistance in mathematics, which could in turn inhibit the potentially positive benefits of such technologies, such as the ability to explore mathematics in broader and more holistic ways as mentioned above.

\subsection{Applying philosophical questions to real-world AI use}

Any content that is a foundational reference for other research carries additional responsibilities, and mathematics is no exception. We can formally certify the validity of any AI-generated mathematical argument; but validity is only one component of value, and there are nuanced value judgments that are necessary in presenting AI-driven research in real-world situations. Which elements of the potentially large body of trivial and non-trivial findings the researcher find particularly interesting and noteworthy to share within and beyond the field of research and how that material is presented to a wider audience has not been standardized among human researchers. There are also uncertainties in how precedence and credit are assigned. AI-assisted research also presents new ethical and legal ramifications and as-yet unanswered questions on the intellectual property rights of AI-generated content (including proofs).

What principles should guide researchers in deciding on the suitability and best application of one AI model or another, or if AI is a good choice at all? In academic domains, it is not unreasonable to make the assumption that most who pursue the path of academic research do so out of a desire to make the world a better place and to make meaningful contributions to it. Mathematicians will want to prioritize those use cases that are the most beneficial to mathematics. Researchers in all fields will frequently want to prioritize not only those uses that benefit their own field but which have cross-disciplinary benefits as well. And it can be given that most who use AI for research purposes at all will want to prioritize uses which benefit humanity over those which harm it. It is consequently important that within the field of AI development, it needs to be highlighted who is benefiting from these tools and what benefits are occurring to help people identify how to responsibly optimize the outcomes as much as possible.

\subsection{Intellectual property and responsibility}

The issue of intellectual property and responsibility (or perhaps, accountability) alone is a minefield and needs careful discussion. When AI is applied to a problem, who is responsible for errors? Who gets credit for insights? These may  not be the same party and may not be parties that are clearly defined.  So far, much of the accumulation of training data for the large language models (LLMs) has been argued (by their developers) as falling under the ``Fair Use'' doctrine. Within the United States, the application of ``fair use'' has some flexibility depending upon (among other things), the purpose of the IP use \cite{harvard}. As a thought experiment, we can consider whether greater benefit merits greater use \cite{weir}. Would it be reasonable to claim it is fair use to draw upon all recorded knowledge in a situation where it is intended to save the world from impending doom? Would such a broad application still apply if it was saving the world from a more distant existential threat (e.g., climate change)? What about if it was ``only'' ending all disease? Or simply eradicating cancer? As all of these are posited beneficial applications of AI, is it then reasonable to grant AI use of all recorded information to make such marvels possible? 

Beyond the problematic argument for an extremely broad interpretation of ``fair use'', clear standards and protocols for the assigning of credit and for citation are desperately needed. AI-use cases will draw not only upon the researcher's data but also the information the AI was previously trained on, the choices of which information the AI was trained on (made by software engineers and designers who may have no interaction with the primary researcher) and, of course, the AI contributions themselves. Is the traditional academic citation system adequate for assigning proper credit in a situation with potentially hundreds or thousands of ``hidden'' contributors, or is it adequate to simply cite the AI model itself?   The undisclosed use of AI to perform a significant portion of the writing in research papers has provoked particularly strong reactions, with many academics viewing such practice as comparable to plagiarism; ironically, this has led some researchers who derive benefit from their tools to conceal their usage from view even further.  It is clear that new professional standards and practices regarding AI disclosure and use will need to be developed\footnote{For an initial discussion of this topic by one of the authors, see \cite{best}.}.

AI is also on the verge of creating potentially widespread circular citation loops, a process humorously dubbed ``citogenesis'' by Randall Munroe\footnote{\url{https://xkcd.com/978}} in 2001.  For instance, following the recent success of AI ``deep research'' tools \cite{openai} in uncovering solutions to open problems that had been buried in obscure literature, the second author helped launch an effort on an mathematical open problem site \cite{bloom} to systematically use these tools to report the known literature on these problems, or the absence of such.  While this added real value to the site, we also found that the deep research tools used these reports as an authoritative source for their search, with the unintended consequence that summarizing these searches on the site interfered with any subsequent use of these tools to turn up genuinely new literature on such problems!  Thus, even in the absence of malicious intent, the increasing power of these tools necessitates a more thorough vetting of the provenance of cited information.

\section{The costs and benefits of AI}

\subsection{Economic and societal impacts: who benefits?}

Given the already significant impact of AI on individuals, as well as its rapid pace of development, it is all too easy  to see a pathway in which AI scales up to present a species-wide existential threat. With any steps forward, developers and other influential individuals need to carefully consider who is benefiting from these advance and who is being harmed by them. We would propose that any further development should prioritize benefit for humanity as a whole and that AI applications should remain directly useful to humans (individually or collectively).

For each individual use-case, assessment should be done to articulate who the intended beneficiaries are. Will this specific AI model or implementation of a model benefit society as a whole or will it only deliver tangible benefits (such as cost savings) to a small group of individuals? AI tools are of such power and complexity that extreme economic gain for a small number of individuals at the cost of millions comes at an intolerable and unacceptable moral cost. We must facilitate implementations of AI which preserve and value the humanity of humans above their commodification.

We need not look far for disastrous outcomes of prioritizing capital over human well-being. Often characterized as arbitrarily anti-technology and anti-progress, the Nottingham textile workers of the early nineteenth century who called themselves ``Luddites'' objected violently to automation displacing them from their jobs and replacing them with lower-skilled and lower-wage workers. The immediate threat to their jobs presented an existential threat to their livelihoods in a harsh economic climate characterized by high unemployment and rampant inflation. While we look back at the automation of the Industrial Revolution as general beneficial to society, these benefits came with real, measurable human costs.

Today, unlike in the Luddites' time, we are already seeing skilled workers replaced not with lower-wage human labor, but with AI. Entry level jobs have  historically been the path to financial and social prosperity for a burgeoning generation of workers. When they simply vanish, opportunity vanishes with them. Despair and resentment build to anger and outrage as humans place themselves in direct opposition to the tools that held promise to improve their quality of life.

As all emerging technologies have some benefits for humanity as a whole, they also come with a real human cost. For a radically disruptive technology like AI, the human costs must be quantified at a local level and a global level and carefully weighed against the benefits. The metrics we use for this assessment are still fuzzy and ill-defined. Do we continue as we have, to look at the monetary gains and losses? Should we be considering the increased access to resources balanced against the resources lost? Do we consider the more intangible benefits of quality of life and happiness, and if so, how do we compare these intangibles against more quantitative gains?

The current business climate unfortunately seeks a \emph{Wunderwaffe} which is being optimized for power and the broadest possible impacts in the hopes that it will be able to outrun any potential problems.  But a failure to quantify the human cost of our emergent technologies does a great disservice to humanity as a whole for the benefit of a select few. Further, the current climate where AI is being implemented simultaneously in virtually every sphere of society, without consideration for whether it provides the end users any meaningful benefit, only serves to alienate and frustrate people in all walks of life. We are already seeing the natural reaction to having a technology imposed on individuals without consent--feeling a loss of control, their first instinct is to reject all AI technologies, even at the risk of throwing the ``baby'' (AI uses that offer quantifiable benefits in their lives) out with the ``bathwater''. If we can instead keep our technology focused on first and foremost, quantifiably improving most or all human lives, we are much less likely to destroy ourselves than if the sole focus of these technologies is the commodification of mechanical labor, digital labor, and human labor.

\subsection{Tallying the costs of AI}

Alongside the direct human costs, no ethical implementation of AI can occur without looking at other more opaque and hidden costs. The most substantial and immediately apparent cost of developing and building an effective AI infrastructure comes from the reality that these technologies, unlike those of the computing revolution of the 1970s, cannot be developed as a hobby or cottage industry--there is no garage of computer parts that a single innovative thinker like Steve Jobs can use to build an empire. The AI models that have been built require a massive investment in hardware, servers, talent, and pre-training long before you can get to a working AI, let alone a profitable one.

A better comparison for the scale that AI requires for development is the transcontinental railroad network built in the US in the latter half of the nineteenth century. The companies that built the railroads had to develop and build a fleet of massive engines, and plan and lay thousands of miles of rail before the first train could quickly and reliably transport goods from Iowa to San Francisco, unlocking the economic returns these companies had gambled on. 

The huge upfront outlay for AI-based technologies has led developers to chase a profit-driven capitalist model, creating a new class of technological elites who wrangle enormous sums of invested capital and managed debt while strategically maneuvering to capture and hold finite resources (in land, energy, water, skilled labor, and such) just as the robber barons of the nineteenth century Gilded Age did. As with that age, the scale of these investments has resulted in massive inequities in economic stability, in access to these technologies, and in general quality of life across the developed world.

Our society has already begun to recognize the significant environmental costs that large-scale AI demands. Heavy energy and water consumption create significant daily challenges for those living in the shadow of the expansive facilities these AI models require. It has been credibly suggested (see, e.g., \cite{cowls}) that AI-generated solutions can be applied to mitigate or eliminate the heavy climate costs of two centuries of human technology use.  And perhaps the marginal costs of operating these tools will decline over time as the infrastructure is built out, and more efficient uses of computation are developed. However, to date, none of the large AI models in operation have provided a solution to even offset their own resource consumption and waste emissions. 

Additionally, it is noteworthy that modern AI tools do not pursue or intuit ``truth'' through manifestation in the physical world, or comprehension of the immutable nature of our reality's physical laws; instead, these models rely heavily on human-generated data, often without attribution, as well as significant amounts of human feedback to iteratively improve itself. Models cannot be built to be less reliant on human intellectual labor without a serious risk of contaminating our collective body of information with AI-generated information. There is a clear limit to how much AI can be used to generate ``new information'' in a domain before AI collapse \cite{shumailov} becomes a serious problem. Without a sufficient amount of genuine content, AI becomes ungrounded from reality, caught up in a mode of thought that has no connection to the real world and significantly hampers the meaningful interactions at the human/AI interface. Mathematics, with its formal verification process, may have a tolerance for higher levels of AI contamination than other domains; but as we have seen, it is not completely immune to this danger.

\subsection{The Digital Divide}

A further significant social cost to consider is the potential for AI technologies to exacerbate existing inequalities or to create new ones. In principle, all humans have the ability to utilize their natural intellectual talents (assuming adequate education and a supportive environment, of course); but the trends in the application of frontier AI models already demonstrate that the large scale AI tools may only be available to well-financed or well-connected research groups, or to individuals who are the most willing to hand over their personal data and look past any ethical concerns regarding the use of such models. This creates a fundamental ``digital divide'' between the AI-have and the AI-have-nots.

It is paramount to prioritize equitable access when AI has the capacity to radically improve research performance, however within the current AI landscape, a second, more nuanced digital divide appears. When the dominant AI models are capitalized, privatized, and competing for finite resources (investment and dependent user base), they are (perhaps unintentionally) incentivized to develop ``spiky'' capabilities to retain a competitive advantage over each other, rather than to provide consistent and even performance in different domains. As individuals are locked into one model over the others due to institutional negotiations and market restraints, we must consider the risk that one model will give a meaningful advantage over another in a particular research sphere, creating divisions even within the subgroup that has reliable and easy access to AI resources.

On the other hand, many of the benefits of AI models in scientific and humanities-based research do not necessarily require the most advanced models. Smaller ``local models'', as well as non-LLM technologies such as proof assistants, demonstrate the capacity to return meaningful results faster and more efficiently than models that necessitate massive data centers processing the sum of all human knowledge. There is significant potential to be able to distill smaller models from the existing larger ones to take advantage of the most advanced AI capabilities with small, user-defined training libraries carefully targeted to the specific area of research interest. Perhaps a diverse array of smaller, more targeted models, maintained by a community of users, could emerge as a viable alternative the current extremely large and expensive models available today. Increased support for such community projects and could help to alleviate the problem of inequitable access. 

While many of these smaller projects would be feasibly developed and run through smaller-scale public and private institutions, industry practitioners and policymakers have called for regulatory actions to create and preserve equitable access to AI technologies \cite{pcast}. As part of that effort there would be significant advantages to investing in the development national or multi-national public-facing coalition for advanced AI research and the development of a large, publicly funded and publicly accessible AI resource (or models) \cite{jones} to readily bring AI access to those individuals and groups who would otherwise be left behind by the private, corporatized models that currently dominate the field. 

\subsection{Harm Reduction}

In the early days of aviation, plane travel was an incredibly unsafe technology, with countless horrific accidents.  Today, it is the safest and most reliable mode of transportation over long distances.  Just as AI has the potential to lead to catastrophic outcomes in the near term, for it to follow a similar trajectory (hopefully with fewer fatal incidents) will require decisive actions to reduce harm. Best practices must be defined \cite{mantegna} and training and regulation designed to enhance the most responsible uses of the technology, while discouraging or banning concealed or harmful ones.

This is a fine needle to thread.  On the one hand, an individual who is using AI assistance cautiously and responsibly might be overtaken in the short term by less scrupulous rivals who are using faster, but more unreliable, AI practices to accelerate their work.  At the same time, such individuals may be derided, condemned, and excluded by cohorts of AI-distrusting peers for even daring to entertain the possibility of incorporating the technology into the workflows of their profession. The current, largely \emph{laissez faire} approach to allowing AI technology to develop at an unchecked pace does not seem promising for such a nuanced, responsible approach to adoption to prevail.

There are some precedents to draw upon for guidance.  The rapid development of Wikipedia in the early twenty-first century initially caused some disruption to educational systems, as many students started blindly incorporating text from that online resource \emph{verbatim} into their assignments, and many instructors reacted by banning the use of that encyclopedic resource. Critiques of Wikipedia's unreliability and potential bias were commonplace.  However, as the site matured, and academia gained familiarity with its strengths and weaknesses, a rough consensus emerged on how to incorporate this resource into education and research.  It became encouraged, or at least condoned, for students and researchers alike to use Wikipedia as a starting point for inquiry on a given topic; and, instead of using its text directly, students are urged to follow up the secondary sources provided by the site, or check them against independent sources of information.  Today, Wikipedia is widely accepted as a useful resource in academia.

Could we reach a similar level of responsible acceptance with AI?  We are cautiously optimistic that this is possible; but it will require sustained effort and clear philosophical guidance.  For instance, we believe it to be a moral and ethical imperative that AI tools should be developed to benefit all (or at least most) humans, rather than a privileged few; that it must create solutions to actual human needs and enhance the quality of life and experience for as many humans as possible; and that the real or potential harms of these tools are recognized, assessed against their benefits, and mitigated whenever possible.  It does not require excessive cynicism to recognize that many of these objectives will not be attained in practice; but debating the system of values we wish these tools to align with is the first step to making it possible to actually achieve these goals. 

As some consensus is (hopefully) found on these values,  then in coordination with the actions above to mitigate the worst impacts of AI, attention must be turned to the greatest source of friction--the interface between AI and humans. To move beyond an uneasy and unstable truce, we need to develop methods to enable individuals to incorporate AI tools into their daily life in ways that feels satisfying and energizing instead of draconian and oppressive. As AI continues to develop and evolve, so to will humanity's uses, interactions, and ultimately relationship with AI need to evolve, from convenient tool to assisting partner to ready collaborator.

\section{The human/AI interface}

\subsection{A short term view: AI as the ``vanilla extract'' of intellectual production}
How should we conceptualize the interface between humanity and AI tools?  In the immediate moment, it is still defensible to view these technologies primarily as curiosities and many users are uncertain as to how to reasonably apply them. 

Our suggested guidance for navigating this current transition is to make a culinary analogy: vanilla extract, a common ingredient in most sweet recipes famous for its nearly universally appealing scent. Ingested by itself, vanilla extract is usually considered extremely unpleasant, but its addition in small amounts is widely regarded as improving and enhancing the other flavors of the dish, even when it cannot be differentiated from them. While it is easy to conclude that more vanilla extract is better, most people who have used it understand there is some upper limit beyond which it ruins the dish entirely \footnote{A notorious thought experiment on Tumblr \cite{vanilla} concluded that a cake that was 44\% vanilla extract would be inedible.}. Most of us do not have a clear sense of what that upper limit actually is, so find it wisest to keep it as a very minor addition. 

Similarly, one could view current AI usage as an optional addition to cognitive workflows: it is interesting to experiment with in moderation--a pass of a human-composed text through an AI language model for suggestions on grammar and phrasing, or a list of bullet points handed to AI to organize into a suggested structure. These light touches, like a small splash of vanilla, will enhance and enrich the character of the work without overwhelming it. AI-content that is utilized to as the core components of a such workflows, however, will not yield desirable, effective, or valuable outcomes. With such a philosophy (and appropriate citation of AI use), there is no immediate need to rethink fundamental assumptions about the role of humans in intellectual pursuits such as mathematics, the sciences, or creative arts.

\subsection{The medium-term: AI on the ``red-team''}
However, as these tools increase in capability and become more broadly adopted, the ability to ``opt out'' of such technologies will diminish.  Even if one personally chooses to actively avoid using AI assistance, colleagues, students, and professional institutions that individual interacts with will increasingly incorporate AI into their own work. Presently, there are serious concerns that entire areas of academic discourse could be drowned out by a flood of low-quality AI-generated content.  In the near term, this can be combated with strict editorial policies to prohibit most forms of AI-generated content; but as these tools become more pervasive and a network of individualized AI agents become more commonplace, a more nuanced approach will become necessary.

In the medium term at least, it will still be possible, and necessary, to devise rules and guidelines to identify the more responsible usages of AI and discourage irresponsible use, without fundamentally changing the humanistic nature of one's field - in short, viewing AI assistance as a tool or junior partner, rather than as a replacement, for human-centered work.   In this case, it can be useful to make a distinction between\footnote{This terminology is inspired by the distinction in cybersecurity between a ``blue team'' that defends a system from attackers, and a ``red team'' that probes for weaknesses.} the ``blue team'' tasks of generating new content and structures, and the ``red team'' task of verifying, testing, or maintaining that content.  AI is relatively safe to utilize in a ``red team'' capacity of reviewing human-generated content for errors or suggested improvements; but with the stochastic unreliability and lack of groundedness of the current and near-term tools, it is unsafe to trust them in any ``blue team'' structural capacity that is beyond the ability of the ``red team'' (which could consist of humans or more automated verification tools, such as formal proof assistants) to verify.  In this philosophy, the emphasis is on managing the potential risks of AI use while still capturing many of its potential benefits,  rather than radically rethinking the fundamental nature of the field.

\subsection{The longer term: is a philosophical retreat inevitable?}
But suppose one looks ahead to a more distant future in which the current weaknesses of AI tools are satisfactorily resolved, and their capabilities now match or exceed that\footnote{This scenario is sometimes referred to as ``Artificial General Intelligence'', although there is no consensus on the precise definition of this term.} of expert humans in all practical dimensions, rendering the risk-management philosophy obsolete.  How will we then respond to the complex philosophical questions raised by the transformative nature of such advanced technologies?

One option is simply to retreat into purely technical frameworks in which these questions are no longer operative.  In mathematics, we have the ``formalist'' viewpoint, where the only objective is to manipulate mathematical symbols according to precise rules.  In the sciences, the pragmatic ``shut up and calculate'' philosophical position plays a similar role; and in the creative arts, one can work as an artisan rather than as an artist, creating works that satisfy the parameters provided by an external client, without passing any judgment on the value of the product.  In each of these cases, no distinction need be made between human-generated work or AI-generated work, so long as the technical specifications of the task are met.

But while technique is certainly an essential component of each of these disciplines, it does not capture the full experience of how mathematics, science, and the arts are conducted in practice, and provides little guidance on such practical questions as how to motivate the next generation of students, or what directions of curiosity-driven research to pursue.  So, one could instead retreat to a radically different position, in which one ascribes an ineffable special status to human intellect or human creativity, permanently distinguishing any activity exercising these human traits such talents at a fundamental level from any artificial replication of that activity, regardless of how accurately the latter could replicate or surpass the former at a technical level.  In this framework, Artificial Intelligence will forever be ``No True Scotsman'': lacking true ``soul'' or ``understanding''.  With the long familiarity with our own species, we are used to humans being unreliable, ``spiky'' in their abilities, and sometimes lucking into successfully achieving a task through random word association and rote memorization; but when AI tools exhibit similar behaviors, one can be inclined to judge them far more harshly, for instance attributing such failings to their inherent nature as ``stochastic parrots''.
But perhaps this position is simply denying an uncomfortable truth: that some portion of our vaunted human capabilities are in fact not that much more sophisticated in nature than the AI algorithms we have now designed to mimic them.  And as AI performance continues to advance, such a human-chauvinistic viewpoint risks degenerating into an increasingly untenable ``god of the gaps'' philosophy, in which an ever-shrinking list of qualities are touted as indicators of essential human achievement that AI is still not yet able to replicate.

A third option, particularly favored by some enthusiasts of these technologies, is to hold that all human cognitive abilities will soon be completely superseded by their AI equivalents, rendering philosophical discussions about the value of human contributions and concerns to mathematics, science and the arts increasingly moot.  In the more extreme versions of this position, the very exercise of human intellect is viewed as an undesirable and tedious activity, which ought to be replaced by automation as quickly as possible, in order to free up time and mental space for more leisurely or hedonistic pursuits.  Obviously, an implementation of this philosophy would carry many risks, such as the degradation of human abilities to the point where our species will become collectively unable to monitor, control, or even understand the actions of that increasingly powerful AIs that we will have delegated our civilization to\footnote{For a vision of what this framework would look like in practice, we suggest the science-fiction film \emph{Wall-E} \cite{wall-e}.}.

There however appears to be some philosophical middle ground between these ``straw-man'' extremes, which can provide useful perspective for emerging paradigms of cooperation and complementary coexistence between humans and AI agents.  One precedent for this can be seen in the world of chess, which was once seen as a quintessential exercise of pure human intellect.  It has now been several decades since any human grandmaster has been able to best a chess engine.  Nevertheless, chess remains a popular and thriving human activity, with chess players incorporating engines into their training, using them to revisit old chess theories and explore new ones, probe for exploits and weaknesses in otherwise invincible AI chess players, or creatively introduce new forms of competition that involve varying levels of AI assistance.  The philosophical questions of what the game of chess actually \emph{is}, and what the value is in playing it, continue to be worth asking; and the currently accepted answers do not closely resemble any of the three extreme positions outlined above. 

\subsection{A Copernican view}
One possibility is to embrace a cognitive analogue of the Copernican revolution in astronomy.  In antiquity, the dominant models of cosmology (insofar as the universe was viewed in mechanistic terms) were geocentric in nature, in which the Earth had a privileged ontological status as the immobile center of the universe, fundamentally distinct in nature from the heavens above or the underworld beneath.    However, multiple advances in astronomy and physics dismantled this view, successively demonstrating over the centuries that the Earth was in fact in motion around its axis, and in orbit around the Sun, with the Sun itself orbiting the center of our galaxy, which in turn was part of an expanding universe that lacks any notion of a spatial center.  Indeed, it became extremely fruitful to adopt a completely opposing philosophical viewpoint, now known as the \emph{Copernican principle}: that the Earth was just one planet among countless others in the universe, receiving no preferential treatment whatsoever from the fundamental laws of nature.   

At first glance, this view feels quite threatening to humanity's emotional attachment to our home planet, but ultimately there is no fundamental contradiction between the universe's disinterest in the planet Earth, and our own strong investment in it; we can still quite justifiability prioritize issues specific to planet Earth over those on other planets, while simultaneously accepting that these other planets exist and would be of comparable importance to their own inhabitants.  Similar revolutions can be seen in the historical development of other sciences, for instance in the Darwinian revolution dethroning the unique status of humans among other constantly evolving species, or the dethroning of the privileged role of Euclidean geometry as a source of synthetic \emph{a priori} truth in mathematics.

Until recently, our species has similarly embraced an intellectual analogue of the geocentric model, in which human intelligence stood at the center of the cognitive universe, thus affording it a special philosophical status.  But now we are discovering (or creating) other ``planets'' of intelligence comparable in many ways to our own, while simultaneously being quite distinct in many aspects.  Instead of denying the existence or importance of these planets, or arguing over which of these planets deserves to be the ``center'',  one can instead accept that both human and artificial intelligences exist in the same ontological category, though with many distinctive differences and complementarities.  While our interests and attachments will still largely be tied to the human intellectual sphere, its relationship with other forms of intelligence can be explored, both for practical purposes of more efficiently achieving various real-world objectives, as well as for more philosophical reasons, such as achieving an external perspective on human cognition that was previously difficult to attain.

\section{Conclusion}
The unstructured, chaotic, and widespread release of AI technology into the world has already dramatically shifted social, intellectual, and economic spheres in ways that are as alarming as they are beneficial. While unquestionably, some kind of collective effort by humanity is needed, whether through regulation, market pressure, or by some as-yet defined force; we have decidedly not yet reached a tipping point from which we cannot extricate ourselves from the high economic and social cost of these new technologies. Approaches to integrating AI into the field of mathematics have just as rapidly demonstrated the promising benefits that AI can bring to academic research, scientific progress, and to humanity at large. The largely objective and verifiable nature of mathematical research presents a unique opportunity to experiment with these new technologies and study the resulting impacts in ways that do not present an ethical or existential risk to the individual or broader society. From the application of AI to mathematics, we are able to explore the pressing philosophical and moral questions of broader global AI use. Further, we can extrapolate potential pathways to relieve the tensions at the AI/human interface and suggest new paradigms of cooperative AI/human thought that respect the unique and valuable qualities that each modality brings to the metaphorical table. Though we will never get the genie back in the bottle, we are optimistic that, as our understandings and action rapidly advance, we can yet clear the smoke away and look toward a bright, if somewhat uncertain, future.

\subsection{Acknowledgments}

We thank Silvia de Toffoli for helpful comments and references.

\bibliographystyle{ieeetr}
\bibliography{AI-Math-Article.bib}


\end{document}